\documentclass[11pt]{article}

\usepackage[margin=1in]{geometry}
\usepackage{enumerate}
\usepackage{amsmath}
\usepackage{amssymb,latexsym}
\usepackage{amsthm}
\usepackage{graphicx}
\usepackage{hyperref}
\usepackage{amscd}

\DeclareMathOperator{\ind}{ind}

\title{Renitent lines}
\author{Bence Csajb\'ok, Peter Sziklai, Zsuzsa Weiner\thanks{The first author was supported by the \'UNKP-19-4 New National Excellence Program of the Ministry for Innovation and Technology and by the National Research, Development and Innovation Office -- NKFIH, grant no. PD 132463. The first and the third author acknowledges the partial support of the National Research, Development and Innovation Office -- NKFIH, grant no. K 124950. The second author is grateful for the partial support of project K 120154 of the National Research, Development and Innovation Fund of Hungary; and for the support of the National Research, Development and Innovation Office within the framework of the Thematic Excellence Program 2021 - National Research Subprogramme: ``Artificial intelligence, large networks, data security: mathematical foundation and applications.'' This work was supported by the Italian National Group for Algebraic and Geometric Structures and their Applications (GNSAGA--INdAM).
}}
\date{}

\newcommand{\cA}{{\mathcal A}}

\newcommand{\cT}{{\mathcal T}}
\newcommand{\cC}{{\mathcal C}}

\newcommand{\cE}{{\mathcal E}}
\newcommand{\cF}{{\mathcal F}}

\newcommand{\cH}{{\mathcal H}}

\newcommand{\cL}{{\mathcal L}}

\newcommand{\cS}{{\mathcal S}}

\newcommand{\F}{{\mathbb F}}
\newcommand{\Fq}{\F_q}

\newcommand{\GF}{\hbox{{\rm GF}}}

\renewcommand{\mod}{\hbox{{\rm mod}\,}}

\newtheorem{theorem}{Theorem}[section]

\newtheorem{lemma}[theorem]{Lemma}
\newtheorem{corollary}[theorem]{Corollary}
\newtheorem{definition}[theorem]{Definition}

\newtheorem{result}[theorem]{Result}
\newtheorem{example}[theorem]{Example}

\newtheorem{remark}[theorem]{Remark}

\DeclareMathOperator{\PG}{{PG}}
\DeclareMathOperator{\AG}{{AG}}

\begin{document}
	\maketitle
	\begin{abstract}
		There are many examples for point sets in finite geometry which behave ``almost regularly'' in some (well-defined) sense, for instance they have ``almost regular'' line-intersection numbers. In this paper we investigate point sets of a desarguesian affine plane, for which there exist some (sometimes: many) parallel classes of lines, such that almost all lines of one parallel class intersect our set in the same number of points (possibly mod $p$, the characteristic). The lines with exceptional intersection numbers are called {\em renitent}, and we prove results on the (regular) behaviour of these renitent lines. As a consequence of our results, we also prove geometric properties of codewords of the $\mathbb{F}_p$-linear code generated by characteristic vectors of lines of $\mathrm{PG}(2,q)$.
	\end{abstract}

	\section{Introduction}
	
	One of the key motivations in the history of finite geometries is the study of symmetric structures, i.e. structures admitting a large symmetry group. These structures (quadrics, Hermitian varieties, subgeometries over a subfield, etc.) are typically very ``regular" when one considers their intersection properties with the subspaces of the ambient geometry; and there exist many ``classification-type" results, stating that an ``intersection-wise very regular" set must be one on the list of the (classical, symmetric) structures. 
	
	Throughout combinatorics, many theorems state that within a certain class of structures, the ``nicest" ones (with respect to some well-defined combinatorial property) are {\em stable}, meaning that if another structure is ``almost as nice", then that must be just a slightly modified version of a nicest one. Such {\em stability results}, as e.g. the Erd\H{o}s--Simonovits theorem (see \cite{Simonovits}) which states the stability of Tur\'an's theorem in extremal graph theory, or the Hilton--Milner (see \cite{Hilton-Milner}) theorem for intersecting families, always describe some deep properties of the class. 
	
	Following this route, it is natural to investigate point sets, which behave ``{\em almost} regularly" with respect to the subspaces of the ambient space. Describing the patterns how little irregularities may occur, reveals some properties of the regular structures.
	
	In this paper we restrict ourselves to point sets of a desarguesian affine plane $\AG(2,q)$, where $q$ is a power of the prime $p$. (The natural but not obvious extensions to other spaces will come in separate papers.)  It may well happen that our point set intersects almost all lines of a parallel class in the same number of points, possibly mod $p$. (E.g. if you take a point set which intersects every line in constant mod $p$ points and delete a few, say $\kappa$ points of it then all the lines, except {\em at most} $\kappa$ per parallel class, still intersect it in constant mod $p$ points.) If it happens for many parallel classes, then one may guess that the reason is that our point set has a hidden structure, i.e. the non-regular intersections may be ``corrected", or at least they also possess some regularity themselves. For instance, in $\AG(2,q)$, $q$ even, an arc of size $q+1$ has almost regular intersection numbers mod 2 in every parallel line class: it meets each line in 0 or 2 points, except one line per each parallel class (which is 1-secant). Now, it is easy to prove that these 1-secant (exceptional) lines are concurrent, which can be interpreted as (i) they are {\em points of a degree $1$ curve in the dual plane}; or (ii) their intersection point can be added to the original point set, which becomes ``regular" this way (i.e. every line meets it in 0 mod 2 points).
	This tiny example illustrates that if you control the irregularities then you may prove that every oval is contained in a hyperoval (i.e. in the ''nice structure"). In general, a $(q+2-k)$-arc is a point set (of size $(q+2-k)$) in $\PG(2,q)$ intersecting every line in at most $2$ points. It is easy to see that such a point set has exactly $k$ $1$-secants through each of its point. Segre's celebrated theorem says that there is a curve of degree $k$ or $2k$, depending on $q$ is even or odd, in the dual plane which contains these $1$-secants. 
	
	There are many papers dealing with problems in view of (ii), for a survey on arcs see \cite{arcsurvey} and for $k \pmod p$ multisets see \cite{kmodp}.
	In this paper we extend and explore this idea. We start with a result, which is from \cite{CsW} and it can be viewed as a generalization of \cite[Theorem 5]{semi}, see also \cite[Proposition 2]{seminuclear} and \cite[Remark 7]{dirszonyi}. We use the usual extension of $\AG(2,q)$ with its {\it line at infinity} $\ell_{\infty}$, containing points called {\it directions}, and the affine lines with slope $d$ all meet at the direction $(d)$, where $d\in\Fq\cup\{\infty\}$.
	
	\begin{lemma}[Lemma of renitent lines  \cite{CsW}]
		\label{renitent}
		Let $\cT$ be a point set of $\AG(2,q)$.
		A line $\ell$ with slope $d$ is called {\em renitent} if there exists an integer $m_d$ such that $|\ell \cap \cT|\not\equiv m_d \pmod p$ but every other line with slope $d$ meets $\cT$ in 
		$m_d$ modulo $p$ points. The renitent lines are concurrent. 
	\end{lemma}
	
	Now we define {\bf renitent lines} in the following, more general setting and then prove various generalizations of the lemma above.
	
	\begin{definition}
		\label{def}
		Let $\cT$ be a multiset of $\AG(2,q)$.
		For some integer $\lambda$,
		%	$\leq (q-1)/2$, 
		a direction $(d)$ is called {\em $(q-\lambda)$-uniform} w.r.t. $\cT$ if there are at least $(q-\lambda)$ affine lines with slope $d$ meeting $\cT$ in the same number of points modulo $p$.
		This number will be called the \emph{typical intersection number} at $(d)$. The rest of the lines with direction $(d)$ will be called {\em renitent} (w.r.t. $\cT$).  
		
		A direction $(d)$ is {\em sharply} $(q-\lambda)$-uniform w.r.t. $\cT$ if it is incident with {\em exactly} $(q-\lambda)$ affine lines meeting $\cT$
		in the same number of points modulo $p$.
	\end{definition}
	
	In the definition above, different directions might have different typical intersection numbers. In most cases we think about $\lambda$ as a small value (hence the name {\em typical}), however in Theorems \ref{regular} and \ref{gcds}, in the prime case we do not need that $\lambda$ is small. When $q$ is non-prime we will usually suppose $\lambda \leq p-1$. In this latter case it is rather automatic that the typical intersection numbers are uniquely determined for each $(q-\lambda)$-uniform direction. 
	
	\begin{remark}
		If $\lambda>(q-1)/2$ then it might happen that for a given $(q-\lambda)$-uniform direction $(d)$ the typical intersection number is not uniquely determined. In this case one typical intersection number has to be fixed w.r.t. the direction $(d)$. 
	\end{remark}

	Under some conditions, if a set of $(q-\lambda)$-uniform directions $\cE_{\lambda}$ is of size at most $q$, then we are able to prove that the renitent lines are contained in an algebraic envelope of relatively small class (i.e. a curve of the dual plane of relatively small degree), see Theorem \ref{regular}. This result resembles Segre's theorem mentioned above.
	
	In Theorem \ref{gcds}, we remove the condition on the size of $\cE_{\lambda}$ and give a more general formulation of Theorem \ref{regular}. 
	The proof relies on $k$-th power sums and the Newton--Girard formulas. 
	If we are more permissive with the renitent lines of $\cT$ then the class of the algebraic envelope might increase. 
	To prove the following result we apply a ``weighted'' version of the Newton--Girard formulas, see Lemma \ref{recursive}.
	
	\medskip
	\noindent
	{\bf Theorem \ref{lambda2}}
	\emph{Take a multiset $\cT$ of $\AG(2,q)$ and an integer $0 < \lambda \leq (q-1)/2$. 
		Let $\cE_{\lambda}$ denote a set of $(q-\lambda)$-uniform directions of size at most $q$. 
		The renitent lines with slope in $\cE_{\lambda}$ are contained in an algebraic envelope of class $\lambda^2$. Furthermore, if a direction is $(q-\lambda)$-uniform, but not sharply uniform, then the line pencil centered at that direction is fully contained in the envelope.}
	\medskip
	
	In Section \ref{sec:4}, we show how to apply the resultant method, cf. \cite{kmodp}, with the help of a polynomial which can detect renitent lines at each uniform direction. 
	For other polynomial techniques used in finite geometry, see for example \cite{survey}. 
	In Section \ref{dualsection}, we dualize some of our results and as a consequence, we also prove geometric properties of codewords of the $\mathbb{F}_p$-linear code generated by characteristic vectors of lines of $\mathrm{PG}(2,q)$. This is a generalization of a result of Blokhuis, Brouwer, Wilbrink \cite[Proposition, p. 66]{BBW}.

	\section{Envelopes of small class}

	In this section our aim is to show that there is an {\it envelope} (i.e. a curve in the dual plane) of relatively small {\it class} (i.e. degree) containing the renitent lines with slope in a given subset $\cF_{\lambda}$ of $(q-\lambda)$-uniform directions.
	There exist at least $\lambda$ common lines of such an envelope and a line pencil centered at a sharply $(q-\lambda)$-uniform direction.
	So if there are $s>0$ sharply $(q-\lambda)$-uniform directions in $\cF_{\lambda}$, then such an envelope has class at least $\min\{\lambda,s\}$ (if it has class less than $\lambda$ then it necessarily contains the $s$ line pencils centered at the sharply $(q-\lambda)$-uniform directions).
	
	Any line set of size at most $|\cF_{\lambda}|\lambda$ is contained in an algebraic envelope of class at most $\lceil\sqrt{2 |\cF_{\lambda}|\lambda}\rceil-1$. 
	Indeed, a $3$-variable homogeneous polynomial $h$ of degree $d$ has $\binom{d+2}{2}$ coefficients. For any point $P\in \PG(2,q)$, the condition that $h$ vanishes at $P$ is equivalent to a linear equation for these coefficients. 
	If $\binom{d+2}{2}$ is larger than the size of a point set $\cS$, then there is a non-trivial solution for the homogeneous system of equations which corresponds to the condition that the points of $\cS$ are zeroes of $h$. 
	In Theorem \ref{regular}, we show that no matter of the size of $\cF_{\lambda}$, with some conditions on the renitent lines we can always construct an envelope of class $\lambda$ containing all renitent lines with slope in $\cF_{\lambda}$.
	The results of this section will rely on the Newton--Girard formulas.

	From now on, we coordinatize the projective plane $\PG(2,q)$ with homogeneous coordinates over $\GF(q)$ in such a way that the line at infinity is $[0:0:1]$, while the $X$- and the $Y$-axes are $[0:1:0]$ and $[1:0:0]$, resp. 
	The direction $(d)$, where $d\in\F_q$, means the point $(1:d:0)$ on the line at infinity. 
	
	Later in this section, we will be able to remove the condition $|\cE_{\lambda}|\leq q$ from Theorem \ref{regular}, see Remark \ref{slightgen}.
	
	\begin{theorem}
		\label{regular}
		Take a multiset $\cT$ of $\AG(2,q)$ and let $\cE_{\lambda}$ denote a set of $(q-\lambda)$-uniform directions of size at most $q$ such that:
		\begin{enumerate}[\rm(i)]
			\item $0<\lambda \leq \min\{q-2,p-1\}$,
			\item for each $(d)\in \cE_{\lambda}$ the renitent lines meet $\cT$ in the same number, say $t_d$, of points modulo $p$,
			\item for each $(d)\in \cE_{\lambda}$ if $m_d$ denotes the typical intersection number at direction $(d)$, then $t_d-m_d$ modulo $p$ does not depend on the choice of $(d)$.
		\end{enumerate}
		Then the renitent lines with direction in $\cE_{\lambda}$ are contained in an algebraic envelope of class $\lambda$.
	\end{theorem}
	
	\begin{proof}
		First we show that the number of renitent lines is the same at each direction $(d)$ of $\cE_{\lambda}$, hence we can find a common value $\lambda' \leq \lambda$ for which all the directions in $\cE_\lambda$ are sharply $(q-\lambda')$-uniform; note that if $\lambda'<\lambda$ then we will find an envelope of even smaller class at the end of the proof.
		Let $(d)$ and $(e)$ denote two directions in $\cE_{\lambda}$ which are sharply $(q-\lambda_d)$-uniform and sharply $(q-\lambda_e)$-uniform, respectively. Then,
		counting points on lines with slope $d$ and lines with slope $e$, we obtain
		\[(q-\lambda_d)m_d+\lambda_d t_d \equiv |\cT| \equiv (q-\lambda_e)m_e+\lambda_e t_e \pmod p,\]
		hence
		\[\lambda_d(t_d-m_d)\equiv \lambda_e(t_e-m_e) \pmod p.\]
		By assumption $(iii)$, $t_d-m_d \equiv t_e - m_e \pmod p$ and $t_d \not\equiv m_d \pmod p$, thus $\lambda_d \equiv \lambda_e \pmod p$. 
		Then $\lambda_d=\lambda_e$ follows from the fact that $0 \leq \lambda_e,\lambda_d \leq \lambda \leq p-1$. 
		
		From now on, we may assume that we had $\lambda=\lambda'$ already, and the directions in $\cE_{\lambda}$ are sharply $(q-\lambda)$-uniform. 
		Since $|\cE_{\lambda}|\leq q$, we may assume $(0:1:0)\notin \cE_{\lambda}$. 
		For each $(1 :d : 0) \in \cE_{\lambda}$ put $(0 : \alpha_1(d) : 1), (0 : \alpha_2(d) : 1),\ldots,(0 : \alpha_{\lambda}(d) : 1)$ for the points of the $Y$-axis on the renitent lines with slope $d$.
		
		Put $s:=|\cT|$ and $\cT=\{(a_i : b_i : 1)\}_{i=1}^s$.
		Next define the polynomials
		\[\pi_k(V):=\sum_{i=1}^s(b_i-a_iV)^k \in \F_q[V]\]
		of degree at most $k$. 
		We project $\cT$ from $(d)$ to the $Y$-axis: the line joining $(1:d:0)$ and $(a_i:b_i:1)$ meets the $Y$-axis at the point $(0:b_i-a_id:1)$, hence for each $(1 : d :0)\in \cE_{\lambda}$ the projected image, i.e. the multiset
		\[M_d:=\{(b_i-a_id)\}_{i=1}^s\]
		contains $m_d\ \mod p$ copies of the values belonging to non-renitent lines and $t_d\ \mod p$ copies of the $\alpha_i(d)$'s, or one may say $m_d\  \mod p$ copies of $\F_q$ and $c\ \mod p$ further copies of $\alpha_i(d)$ for $1 \leq i \leq \lambda$, where $c\in \{1,\ldots,p-1\}$ is an integer such that $c \equiv t_d-m_d \pmod p$. %c_d \neq 0 
		Since $\sum_{\gamma\in \F_q} \gamma^k=0$ for $0\leq k \leq q-2$ and since $\pi_k(d)$ is the $k$-th power sum of the elements in $M_d$, for $0\leq k \leq q-2$ it holds that for 
		$(1:d:0)\in \cE_{\lambda}$
		\begin{equation}
			\label{powermindig2}
			\pi_k(d)=c\sum_{i=1}^{\lambda} \alpha_i(d)^k.
		\end{equation}
		Denote by $\sigma_i(X_1,\ldots,X_{\lambda})$ the $i$-th elementary symmetric polynomial in the variables $X_1,\ldots, X_{\lambda}$.
		Also, for any integer $i\geq 0$ and $d\in \cE_{\lambda}$ put
		\[\sigma_i(d)=\sigma_i(\alpha_1(d),\ldots,\alpha_{\lambda}(d)).\]
		For $p-1 \geq j \geq 1$ define the following polynomial of degree at most $j$:
		\[S_j(V):=(-1)^j\sum_{\substack{n_1+2n_2+\ldots+jn_j=j\\ n_1,n_2,\ldots,n_j\geq 0}}\ \ \prod_{i=1}^j\frac{(-\pi_i(V)/c)^{n_i}}{n_i!i^{n_i}}\in \F_q[V].\]
		Then for $\min \{q-2, p-1\}\geq j\geq 1$ from \eqref{powermindig2} and from the Newton-Girard identities it follows that 
		$S_j(d)=\sigma_j(d)$ for each $(1:d:0)\in \cE_{\lambda}$.
		
		Consider the affine curve of degree $\lambda$ defined by
		\[f(U,V):=U^\lambda-S_1(V)U^{\lambda-1}+S_2(V) U^{\lambda-2}-\ldots+(-1)^{\lambda-1}S_{\lambda-1}(V)U+(-1)^\lambda S_{\lambda}(V).\]
		Then the projective curve of degree $\lambda$ defined by the equation  $g(U,V,W):=W^\lambda f(U/W,V/W)$ contains the point $(\alpha_i(d) : d : 1)$ for each $d\in \cE_{\lambda}$ and $1 \leq i \leq \lambda$.
		Indeed, 
		\[g(U,d,1)=U^\lambda-\sigma_1(d)U^{\lambda-1}+\sigma_2(d) U^{\lambda-2}-\ldots+(-1)^{\lambda-1}\sigma_{\lambda-1}(d)U+(-1)^\lambda\sigma_\lambda(d)=\]
		\[\prod_{i=1}^{\lambda}(U-\alpha_i(d)).\]
		It follows that the lines $[d:-1:\alpha_i(d)]$ are contained in an algebraic envelope of class $\lambda$.
	\end{proof}
	
	Let $U$ be an affine $(q-n)$-set and denote by $D_U$ the set of directions determined by $U$ (that is, directions incident with at least one line meeting $U$ in at least $2$ points). Sz\H onyi proved in \cite[Theorem 4]{dirszonyi} that the affine lines which do not meet $U$ and are incident with a not determined direction are contained in an algebraic envelope of class $n$. If $\lambda \leq \min\{q-2,p-1\}$, then this is a special case of our Theorem \ref{regular}, with $\lambda=n$, $\cE_n=\ell_{\infty} \setminus D_U$ (for $d\in \cE_n$ it holds that $t_d=0$ and $m_d=1$). 
	To construct the envelope, Sz\H onyi's proof relies on the fact that a typical line intersects $U$ in $1$ point exactly, while a renitent line in $0$ point. Our theorem considers a more general setting, but the cost of it is that we need an upper bound on $\lambda$.
	
	As in \cite{dirszonyi}, if $\lambda \leq \sqrt{q}/2$ and $|\cE_{\lambda}|\geq (q+1)/2$, then \cite[Proposition 2]{dirszonyi} yields that the envelope constructed in Theorem \ref{regular} is the product of $n$ pencils. (For an improved version of \cite[Proposition 2]{dirszonyi} see \cite[Lemma 3.2]{dirsziklai}.)   This actually happens as the next example shows.
	
	\begin{example}
		\label{triviex}
		In $\AG(2,q)$ consider an $m$ mod $p$ multiset and remove $\lambda$ of its points. 
		These $\lambda$ points determine at most $\binom{\lambda}{2}$ directions which we will denote by $D$.
		Then the points of $\ell_{\infty}\setminus D$ are $(q-\lambda)$-uniform with typical intersection number $m$ and the renitent lines (lines meeting $U$ in $m-1$ mod $p$ points)
		are contained in the product of the $\lambda$ pencils centered at the removed points.
	\end{example}

	As the next example shows, Theorem \ref{regular} does not hold when $\lambda = p^2-p$, hence we certainly need some restriction on the value $\lambda$. In Section \ref{dualsection}, we show that the bound $\lambda \leq p-1$ is necessary when $q=p^2$.
	
	\begin{example}
		In $\AG(2,q)$, $q=p^h$, $p$ prime, $h>2$, consider a subplane $\PG(2,p)$ and delete $p+1$ collinear points from it. Let $\cT$ denote the point set obtained, $\cT=\{Q_1,Q_2,\ldots,Q_{p^2}\}$. Denote by $P_1,\ldots,P_{p^2+p}$ the points of $\ell_{\infty}$ through which there passes a $p$-secant of $\cT$.
		These points are incident with $p^2-p$ affine lines meeting $\cT$ in $1$ point and with $q-p^2+p$ affine lines meeting $\cT$ in $0$ modulo $p$ points. 
		
		We claim that the $(p^2+p)(p^2-p)$ renitent lines are not contained in an envelope of class $p^2-p$. 
		Indeed, the point $Q_i$ is incident with $p^2+p-(p+1)$ lines containing a point $P_j$ for some $j \in [1, p^2+p]$ and intersecting $\cT$ in $1$ modulo $p$ points, i.e with this many renitent lines; so the pencil with center $Q_i$ would be a component of such an envelope for each $i$. A contradiction by degree considerations.    
	\end{example}
	
	The next theorem seems to have rather artificial conditions and settings. But this provides a middle ground between the kind of best bound $\deg\leq\lambda$ of Theorem \ref{regular} with strong assumptions and the weaker bound
	$\deg\leq\lambda^2$ of Theorem \ref{lambda2} holding in general.
	
	\begin{theorem}
		\label{gcds}
		Take a multiset $\cT$ of $\AG(2,q)$ and let $\cF_{\lambda}$ denote a set of $(q-\lambda)$-uniform directions.
		For each $(d)\in \cF_{\lambda}$ denote the typical intersection number by $m_d$ and denote the intersection numbers of the renitent lines by
		$t_{d,1},t_{d,2},\ldots,t_{d,\lambda_d}$, for some $0<\lambda_d\leq \lambda$. 
		For $c\in \F_p\setminus \{0\}$ define the integers $\lambda_{d,i}(c)\in \{1,\ldots,p-1\}$ such that 
		\[c \lambda_{d,i}(c)\equiv t_{d,i} - m_d \pmod p\]
		and assume that
		\begin{equation}
			\label{bonyi}
			\Lambda_d(c):= \sum_{i=1}^{\lambda_d}\lambda_{d,i}(c)\leq \min\{q-2,p-1\}
		\end{equation}
		holds for each $(d)\in \cF_{\lambda}$ (note that the sum is taken over natural numbers).  
		Then $\Lambda(c):=\Lambda_d(c)$ does not depend on $d$ and the renitent lines with direction in $\cF_{\lambda}$ are contained in an algebraic envelope of class $\Lambda(c)$. 
		If $\ell$ is a renitent line with slope $d$ and with intersection number $t_{d,i}$, then the intersection multiplicity of the pencil centered at $(d)$ with $\cF_{\lambda}$ at $\ell$ is $\lambda_{d,i}(c)$.
		%Moreover, the intersection multiplicity of this envelope with the pencil centered at $(d)$ at a renitent line incident with $(d)$ and with intersection number $t_{d,i}$ is $\lambda_{d,i}(c)$. 
	\end{theorem}
	
	\begin{proof}
		Throughout the proof, we will fix $c$ and so neglect it from $\lambda_{d,i}(c)$, $\Lambda_d(c)$ and $\Lambda(c)$.
		As before, first we show that $\Lambda_d$ does not depend on $(d)$ of $\cF_{\lambda}$. 
		Let $(e_1)$ and $(e_2)$ denote two directions in $\cF_{\lambda}$. 
		Then
		\[q m_{e_1}+\Lambda_{e_1} c \equiv |\cT| \equiv q m_{e_2}+\Lambda_{e_2} c \pmod p,\]
		hence
		\[\Lambda_{e_1}\equiv \Lambda_{e_2} \pmod p.\]
		Then $\Lambda_{e_1}=\Lambda_{e_2}$ follows from the fact that $0 < \Lambda_{e_1},\Lambda_{e_2} \leq p-1$. 
		
		First we prove the assertion for any $\cE_{\lambda} \subseteq \cF_{\lambda}$ such that $|\cE_{\lambda}|\leq q$. In this case, we may assume $(0:1:0)\notin \cE_{\lambda}$. 
		For each $(1 :d : 0) \in \cE_{\lambda}$ put $(0 : \alpha_1(d) : 1), (0 : \alpha_2(d) : 1),\ldots,(0 : \alpha_{\lambda_d}(d) : 1)$ for the points of the $Y$-axis on the renitent lines with slope $d$.
		
		Put $s:=|\cT|$ and $\cT=\{(a_i : b_i : 1)\}_{i=1}^s$.
		Next define the polynomials
		\[\pi_k(V):=\sum_{i=1}^s(b_i-a_iV)^k \in \F_q[V]\]
		of degree at most $k$. For each $(1 : d :0)\in \cE_{\lambda}$ the multiset
		\[M_d:=\{(b_i-a_id)\}_{i=1}^s\]
		contains $m_d$ modulo $p$ copies of $\F_q$ and $t_{d,i}-m_d$ modulo $p$ further copies of $\alpha_i(d)$ for $1 \leq i \leq \lambda_d$. %c_d \neq 0 
		Since $\sum_{\gamma\in \F_q} \gamma^k=0$ for $0\leq k \leq q-2$ and since $\pi_k(d)$ is the $k$-th power sum of the elements in $M_d$, for $0\leq k \leq q-2$ it holds that for 
		$(1:d:0)\in \cE_{\lambda}$
		\begin{equation}
			\label{powermindig3}
			\pi_k(d)=c\sum_{i=1}^{\lambda_d} \sum_{j=1}^{\lambda_{d,i}} \alpha_i(d)^k.
		\end{equation}
		Denote by $\sigma_i(X_1,\ldots,X_{\Lambda})$ the $i$-th elementary symmetric polynomial in the variables $X_1,\ldots, X_{\Lambda}$.
		Also, for any integer $i\geq 0$ and $d\in \cE_{\lambda}$ put
		\[\sigma_i(d)=\sigma_i({\underbrace{\alpha_1(d),\ldots,\alpha_1(d)}_\text{$\lambda_{d,1}$ times}},\ldots,{\underbrace{\alpha_j(d),\ldots,\alpha_j(d)}_\text{$\lambda_{d,j}$ times}},\ldots,{\underbrace{\alpha_{\lambda_d}(d),\ldots,\alpha_{\lambda_d}(d)}_\text{$\lambda_{d,\lambda_d}$ times}}).\]
		For $p-1 \geq j \geq 1$ define the following polynomial of degree at most $j$:
		\[S_j(V):=(-1)^j\sum_{\substack{n_1+2n_2+\ldots+jn_j=j\\ n_1,n_2,\ldots,n_j \geq 0}}\ \ \prod_{i=1}^j\frac{(-\pi_i(V)/c)^{n_i}}{n_i!i^{n_i}}\in \F_q[V].\]
		Then for $\min \{q-2, p-1\}\geq j\geq 1$ from \eqref{powermindig3} and from the Newton--Girard identities it follows that 
		$S_j(d)=\sigma_j(d)$ for each $(1:d:0)\in \cE_{\lambda}$.
		
		Consider the affine curve of degree $\Lambda$ defined by
		\[f(U,V):=U^{\Lambda}-S_1(V)U^{\Lambda-1}+S_2(V) U^{\Lambda-2}-\ldots+(-1)^{\Lambda-1}S_{\Lambda-1}(V)U+(-1)^\Lambda S_{\Lambda}(V).\]
		Then the projective curve of degree $\Lambda$ defined by the equation  $g(U,V,W):=W^\Lambda f(U/W,V/W)$ contains
		the point $(\alpha_i(d) : d : 1)$ for each $d\in \cE_{\lambda}$ and $1 \leq i \leq \lambda_d$, with multiplicity $\lambda_{d,i}$. Indeed, 
		\[g(U,d,1)=U^\Lambda-\sigma_1(d)U^{\Lambda-1}+\sigma_2(d) U^{\Lambda-2}-\ldots+(-1)^{\Lambda-1}\sigma_{\Lambda-1}(d)U+(-1)^\Lambda\sigma_\Lambda(d)=\]
		\[\prod_{i=1}^{\lambda_d}(U-\alpha_i(d))^{\lambda_{d,i}}.\]
		It follows that the point $(\alpha_i(d):d:1)$ lies on the curve defined by $g$. 
		Note that $(V-d)$ cannot divide $f(U,V)$ and hence the ``horizontal'' line $[0:1:-d]$ cannot be a component of the curve defined by $f$. 
		The intersection multiplicity of the curve defined by $f$ and $[0:1:-d]$ is $\lambda_{d,i}$ at the point $(\alpha_i(d):d:1)$. 
		It follows that the lines $[d:-1:\alpha_i(d)]$ are contained in an algebraic envelope of class $\Lambda$ and the intersection multiplicity of this envelope with the pencil centered at $(1:d:0)$ is $\lambda_{d,i}$ at the line $[d:-1:\alpha_i(d)]$.
		
		Now assume that $|\cF_{\lambda}|=q+1$.
		Apply the argument above for two distinct subsets of $\cF_{\lambda}$, both of them of size $q$. Denote them by $\cE_{\lambda}$ and $\cE'_{\lambda}$ and 
		denote the corresponding curves (in the dual plane) of degree $\Lambda$ by $\cC$ and $\cC'$, respectively. Our aim is to prove that these two curves coincide. 
		Put $\cC=\cH\cdot \cA$ and $\cC'=\cH\cdot \cA'$, where $\cA$ and $\cA'$ do not have a common component. Denote the degree of $\cH$ by $h$ and suppose to the contrary 
		that $h<\Lambda$. 
		
		For $(d)\in \cE_{\lambda} \cap \cE'_{\lambda}$ denote the line $[0:1:-d]$ by $\ell_d$ and the point $(\alpha_i(d):d:1)$ by $P_{d,i}$. Recall that $\ell_d$ cannot be a component of $\cC$ or $\cC'$. 
		The intersection multiplicity $\mathrm{I}(\cC\cap \ell_d,P_{d,i})=\lambda_{d,i}$ equals 
		$\mathrm{I}(\cH\cap \ell_d,P_{d,i})+\mathrm{I}(\cA\cap \ell_d,P_{d,i})$, thus
		\[\mathrm{I}(\cA\cap \ell_d,P_{d,i})=\lambda_{d_i}-\mathrm{I}(\cH\cap \ell_d,P_{d,i}).\]
		Replacing $\cA$ by $\cA'$, the same argument yields that
		\[\mathrm{I}(\cA'\cap \ell_d,P_{d,i})=\lambda_{d_i}-\mathrm{I}(\cH\cap \ell_d,P_{d,i}).\]
		By \cite[Lemma 9.2 and p.\ 87]{Seidenberg} or \cite[Lemma 10.4]{polybook} it follows that $\mathrm{I}(\cA\cap \cA',P_{d,i})\geq \lambda_{d_i}-\mathrm{I}(\cH\cap \ell_d,P_{d,i})$.
		Then 
		\begin{equation}
			\label{lotlenseg}
			\sum_{i=1}^{\lambda_d} \mathrm{I}(\cA \cap \cA',P_{d,i}) \geq \sum_{i=1}^{\lambda_d} \lambda_{d_i}- \sum_{i=1}^{\lambda_d}\mathrm{I}(\cH\cap \ell_d,P_{d,i})=\Lambda-\sum_{i=1}^{\lambda_d}\mathrm{I}(\cH\cap \ell_d,P_{d,i})\geq \Lambda-h,
		\end{equation}
		since $\sum_{i=1}^{\lambda_d}\mathrm{I}(\cH\cap \ell_d,P_{d,i}) \leq \deg \cH=h$. The inequality \eqref{lotlenseg} holds for each $(d)\in \cE_{\lambda} \cap \cE'_{\lambda}$ and hence
		\[\sum_{P \in \cA \cap \cA'} \mathrm{I}(\cA \cap \cA',P) \geq \sum_{(d)\in \cE_{\lambda} \cap \cE'_{\lambda}} \sum_{i=1}^{\lambda_d} \mathrm{I}(\cA\cap \cA',P_{d,i})\geq (q-1)(\Lambda-h).\]
		On the other hand, By B\'ezout's theorem $\sum_{P \in \cA \cap \cA'} \mathrm{I}(\cA \cap \cA',P) \leq \deg \cA \cdot \deg \cA' = (\Lambda-h)(\Lambda-h)$.
		
		This contradiction proves $h=\Lambda$, that is, $\cC=\cC'$.
	\end{proof}
	
	\begin{remark}
		\label{slightgen}
		If in Theorem \ref{gcds} we assume $t_{d,1}=t_{d,2}=\ldots=t_{d,\lambda_d}=:t_d$ for each $(d)\in \cF_{\lambda}$, we also assume that 
		$t_d-m_d$ does not depend on the choice of $d$, and further assume $0<\lambda\leq \min\{q-2,p-1\}$, then with the choice $c \equiv t_d-m_d \pmod p$ we obtain Theorem \ref{regular} without the restriction 
		$|\cE_{\lambda}|\leq q$. 
	\end{remark}
	
	\begin{remark}
		If $|\cT|\equiv 0 \pmod p$ then Theorem \ref{gcds} cannot be applied. Indeed, in that case for each $(d)\in \cF_{\lambda}$, 
		\[\sum_{i=1}^{\lambda_d} (t_{d,i} - m_d)\equiv 0 \pmod p\]
		and hence $\Lambda_d=\sum_{i=1}^{\lambda_d} \lambda_{d,i} \equiv 0 \pmod p$, which is not possible if $\Lambda_d \leq p-1$ (by definition $\lambda_{d,i}>0$ for each $i$).
	\end{remark}
	
	\begin{remark}
		For $c=1,2,\ldots,p-1$, the values of $\Lambda_d(c)$ give different residues modulo $p$ and this residue is the same for every choice of $d$. 
		However, it is possible that for some of the directions, $\Lambda_d(c) \leq \min\{q-2,p-1\}$ does not hold.
		
		Indeed, put for example $p=5$, $\lambda=2$ and assume that the typical intersection number is $0$ at each direction. 
		Also, assume that both of the two renitent lines with direction $(d_1)\in \cF_{\lambda}$ meet $\cT$ in $1$ modulo $p$ points and the two renitent lines with direction $(d_2)\in \cF_{\lambda}$ meet $\cT$ in $3$ and in $4$ points modulo $p$.
		Then, with $c=1$ we obtain $\Lambda_{d_1}(c)=2$ and $\Lambda_{d_2}(c)=7$, so the result cannot be applied. On the other hand, 
		with $c=3$ we obtain $\Lambda_{d_1}(c)=\Lambda_{d_2}(c)=4$ and hence the result can be applied if $q>5$. 
		This also shows that sometimes it might be convenient to choose $\cF_{\lambda}$ not as the set of all $(q-\lambda)$-uniform directions, but as a subset of them. 
	\end{remark}
	
	\begin{example}
		In Theorem \ref{gcds}, put $\lambda=3$ and assume that $m_d=1$ and the renitent lines meet $\cT$ modulo $p$ in the multiset $\{3,3,5\}$ for each $(d)\in \cF_{\lambda}$. 
		Note that this implies $p\neq 2$. 
		With the choice $c=2$, it follows that the renitent lines are contained in a curve of degree $\Lambda=(3-1)/2+(3-1)/2+(5-1)/2=4$ whenever $4\leq \min\{q-2,p-1\}$.
	\end{example}
	
	One might think that to obtain a curve of the lowest degree, the best option is to chose $c$ as the greatest common divisor of the values $t_{d,i}-m_d$.
	The next example disproves this belief. 
	
	\begin{example}
		In Theorem \ref{gcds},  put $\lambda=2$, $p=13$, and assume that $m_d=1$ and the renitent lines meet $\cT$ modulo $13$ in the set 
		$\{2,8\}$ for each $(d)\in \cF_{\lambda}$. 
		With the choice $c=1$, it follows that the renitent lines are contained in a curve of degree $\Lambda=(2-1)+(8-1)=8$.
		With the choice $c=7$, it follows that the renitent lines are contained in a curve of degree $\Lambda=1/7+7/7=2+1=3$.
	\end{example}

	\section{The general case}
	\label{sec:3}
	
	First we prove a recursion connecting elementary symmetric polynomials with ``weighted'' power sums.
	
	\begin{lemma}
		\label{recursive}
		Let $\sigma_k=\sigma_k(X_1,\ldots,X_\lambda)$ denote the $k$-th elementary symmetric polynomial in the variables $X_1,\ldots,X_{\lambda}$. For some field elements $c_1,c_2,\ldots,c_{\lambda}$ put
		\[P_k=P_k(X_1,\ldots,X_\lambda)=\sum_{i=1}^{\lambda} c_i X_i^k.\]
		Then for any integer $j\geq 0$ it holds that
		\[P_{\lambda+j}=P_{\lambda+j-1}\sigma_1-P_{\lambda+j-2}\sigma_2+\ldots+(-1)^{\lambda+1}P_j \sigma_{\lambda}.\]
	\end{lemma}
	\begin{proof}
		Note that
		\[Y^j\prod_{i=1}^\lambda(Y-X_i)=
		Y^{\lambda+j}-Y^{\lambda+j-1}\sigma_1+Y^{\lambda+j-2}\sigma_2-\ldots+(-1)^{\lambda}Y^j\sigma_\lambda.\]
		It follows that 
		\[c_1(X_1^{\lambda+j}-X_1^{\lambda+j-1}\sigma_1+X_1^{\lambda+j-2}\sigma_2-\ldots+(-1)^{\lambda}X_1^j\sigma_\lambda)=0,\]
		\[c_2(X_2^{\lambda+j}-X_2^{\lambda+j-1}\sigma_1+X_2^{\lambda+j-2}\sigma_2-\ldots+(-1)^{\lambda}X_2^j\sigma_\lambda)=0,\]
		\[\vdots\]
		\[c_\lambda(X_\lambda^{\lambda+j}-X_\lambda^{\lambda+j-1}\sigma_1+X_\lambda^{\lambda+j-2}\sigma_2-\ldots+(-1)^{\lambda}X_\lambda^j\sigma_\lambda)=0.\]
		Summing up both sides above yields the assertion.
	\end{proof}
	
	\begin{lemma}
		\label{matrixos}
		As before, put
		\[P_k=P_k(X_1,\ldots,X_{\lambda})=\sum_{i=1}^{\lambda} c_i X_i^k,\]
		and define the $\lambda\times \lambda$ matrix
		\[H=H(X_1,\ldots,X_\lambda)=\begin{pmatrix}
			P_{\lambda-1} & P_{\lambda-2} & \ldots & P_0 \\
			P_{\lambda} & P_{\lambda-1} & \ldots & P_1 \\
			\vdots & \vdots & \vdots & \vdots \\
			P_{2\lambda-2}& P_{2\lambda-3}& \ldots & P_{\lambda-1} 
		\end{pmatrix}.\]
		Then $\det H = (-1)^{\lambda(\lambda-1)/2}c_1 c_2 \ldots c_{\lambda} \prod_{1\leq i < j \leq \lambda} (X_i-X_j)^2$.
	\end{lemma}
	
	\begin{proof}
		It follows from the properties of Vandermonde matrices and from the fact that
		\[H=\begin{pmatrix}
			c_1 & c_2 &\ldots & c_\lambda \\
			c_1 X_1 & c_2 X_2 & \ldots & c_\lambda X_\lambda \\
			\vdots \\
			c_1 X_1^{\lambda-1} & c_2 X_2^{\lambda-1} & \ldots & c_\lambda X_\lambda^{\lambda-1}
		\end{pmatrix}
		\begin{pmatrix}
			X_1^{\lambda-1} & X_1^{\lambda-2} &\ldots & 1 \\
			X_2^{\lambda-1} & X_2^{\lambda-2} &\ldots & 1 \\
			\vdots \\
			X_\lambda^{\lambda-1} & X_\lambda^{\lambda-2} &\ldots & 1 \\
		\end{pmatrix}.\]
		
	\end{proof}
	
	Now we are ready to prove an upper bound for the degree of the curve, under very general conditions.
	
	\begin{theorem}
		\label{lambda2}
		Take a multiset $\cT$ of $\AG(2,q)$ and an integer $0 < \lambda \leq (q-1)/2$. 
		Let $\cE_{\lambda}$ denote a set of $(q-\lambda)$-uniform directions of size at most $q$ and assume that there is at least one sharply $(q-\lambda)$-uniform direction in $\cE_{\lambda}$. 
		Then the renitent lines with slope in $\cE_{\lambda}$ are contained in an algebraic envelope of class $\lambda^2$ with the following properties:
		\begin{enumerate}[(1)]
			\item if a direction $(u)$ of $\cE_{\lambda}$ is $(q-\lambda_u)$-uniform for some $0 \leq \lambda_u<\lambda$ then the line pencil centered at $(u)$ is fully contained with multiplicity $(\lambda-\lambda_u)$ in the envelope;
			\item pencils centered at sharply $(q-\lambda)$-uniform directions are not contained in the envelope.
		\end{enumerate} 
	\end{theorem}
	\begin{proof}
		For each $(1 :d : 0) \in \cE_{\lambda}$ denote by $\lambda_d$ the number of renitent lines with slope $d$ and denote by $(0 : \alpha_1(d) : 1), (0 : \alpha_2(d) : 1),\ldots,(0 : \alpha_{\lambda_d}(d) : 1)$ the points of the $Y$-axis on these lines (first we assume $\lambda_d > 0$, the $\lambda_d=0$ case is very simple and it is treated in the last paragraph of the proof). 
		Also, denote the typical intersection number at $(d)$ by $m_d$. 
		Put $s:=|\cT|$ and $\cT=\{(a_i : b_i : 1)\}_{i=1}^s$.
		Next define the polynomials
		\[\pi_k(V):=\sum_{i=1}^s(b_i-a_iV)^k \in \F_q[V]\]
		of degree at most $k$. As in the earlier proofs, for any $(1 : d :0)\in \cE_{\lambda}$ the multiset
		\[M_d:=\{(b_i-a_id)\}_{i=1}^s\]
		contains $m_d$ modulo $p$ copies of $\F_q$ and $c_i(d)\not\equiv 0$ modulo $p$ further copies of $\alpha_i(d)$ for $1 \leq i \leq \lambda_d$. %c_d \neq 0 
		Since $\sum_{g\in \F_q} g^k=0$ for $0\leq k \leq q-2$ and since $\pi_k(d)$ is the $k$-th power sum of $M_d$, for $0\leq k \leq q-2$ it holds that
		\begin{equation}
			\label{powermindig}
			\pi_k(d)=\sum_{i=1}^{\lambda_d} c_i(d)\alpha_i(d)^k.
		\end{equation}
		Note that $\pi_k(d)$ is as $P_k(\alpha_1(d),\alpha_2(d),\ldots,\alpha_{\lambda_d}(d))$ (with $c_i=c_i(d)$) in Lemma \ref{recursive}. 
		For any integer $i\geq 0$ and $(1:d:0)\in \cE_{\lambda}$ put
		\[\sigma_i(d)=\sigma_i(\alpha_1(d),\ldots,\alpha_{\lambda_d}(d)),\]
		where for $i>\lambda_d$ we define $\sigma_i(d)$ to be $0$.
		
		For $j\geq 0$, $q-2 \geq \lambda_d+j$ and $(1:d:0)\in \cE_{\lambda}$ Lemma \ref{recursive} yields
		\begin{equation}
			\label{main}
			\pi_{\lambda_d+j}(d)=\pi_{\lambda_d+j-1}(d) \sigma_1(d)-\pi_{\lambda_d+j-2}(d)\sigma_2(d)+\ldots+(-1)^{\lambda_d+1}\pi_j(d)\sigma_{\lambda_d}(d).
		\end{equation}
		Define
		\[	H(V)=\begin{pmatrix}
			\pi_{\lambda-1}(V) & \pi_{\lambda-2}(V) & \ldots &\pi_0(V) \\
			\pi_{\lambda}(V) & \pi_{\lambda-1}(V) & \ldots & \pi_1(V) \\
			\vdots & \vdots & \vdots & \vdots \\
			\pi_{2\lambda-2}(V) & \pi_{2\lambda-3}(V)& \ldots & \pi_{\lambda-1}(V) 
		\end{pmatrix}.\]
		Then for each $(1:d:0)\in \cE_{\lambda}$, since $2\lambda-1 \leq q-2$, by \eqref{main} applied for $\lambda-\lambda_d \leq j \leq 2\lambda-\lambda_d-1$, we obtain
		\begin{equation}
			\label{mainmatrix}
			H(d)\begin{pmatrix}
				\sigma_1(d)\\
				-\sigma_2(d)\\
				\vdots\\
				(-1)^{\lambda+1}\sigma_\lambda(d)
			\end{pmatrix}	=
			\begin{pmatrix}
				\pi_\lambda(d)\\
				\pi_{\lambda+1}(d)\\
				\vdots\\
				\pi_{2\lambda-1}(d)
			\end{pmatrix}.
		\end{equation}
		Denote by $H_i(V)$ the matrix obtained from $H(V)$ by replacing its $i$-th column with 
		\[C(V):=(\pi_{\lambda}(V),\pi_{\lambda+1}(V),\ldots,\pi_{2\lambda-1}(V))^T.\] 
		
		Put $S(V)=\det H(V)$ and $S_i(V)=\det H_i(V)$. Then $\deg S \leq \lambda(\lambda-1)$ and $\deg S_i(V)\leq \lambda^2-\lambda+i$. Also, $S(d)=\det H(d)$ and $S_i(d)=\det H_i(d)$ for $(1:d:0)\in \cE_{\lambda}$.
		Consider the affine curve $\cC$ of degree at most $\lambda^2$ defined by the equation 
		\begin{equation} \label{revX}
			f(U,V):=
			S(V)U^\lambda-S_1(V)U^{\lambda-1}-S_2(V) U^{\lambda-2}-\ldots-S_{\lambda-1}(V)U- S_{\lambda}(V).
		\end{equation}
		
		By Lemma \ref{matrixos}, if $(1:d:0)\in \cE_{\lambda}$ is sharply $(q-\lambda)$-uniform, then the determinant of $H(d)$ is 
		$(-1)^{\lambda(\lambda-1)/2}\prod_{i=1}^{\lambda} c_i(d) \prod_{1 \leq i < j \leq \lambda} (\alpha_i(d)-\alpha_j(d))^2$
		and hence $H(d)$ is invertible and the existence of a sharply $(q-\lambda)$-uniform direction ensures that $S(V)$ is not the zero-polynomial. 
		Then by Cramer's rule $(-1)^{i+1}\sigma_i(d) = \det H_i(d) / \det H(d)=S_i(d)/S(d)$ for $1\leq i \leq \lambda$.
		
		Hence if $(1:d:0) \in \cE_{\lambda}$ is sharply $(q-\lambda)$-uniform (that is $\lambda_d = \lambda$), then $(\alpha_i(d),d)$ is a point of $\cC$ defined in \eqref{revX} for each $1\leq i \leq \lambda_d$,
		since in this case $f(U,d)$ equals
		\[S(d)(U^{\lambda}-\sigma_1(d)U^{\lambda-1}+\sigma_2(d) U^{\lambda-2}-\ldots+(-1)^{\lambda-1}\sigma_{\lambda-1}(d)U+(-1)^{\lambda}\sigma_{\lambda}(d))=\]
		\[S(d)\prod_{i=1}^{\lambda}(U-\alpha_i(d)).\]
		Now consider $(1:d:0)\in \cE_{\lambda}$ such that $\lambda_d < \lambda$. 
		To show that the pencil with carrier $(d)$ is contained with multiplicity $(\lambda-\lambda_d)$ in the envelope, it is enough to prove that $(V-d)^{\lambda-\lambda_d}$ divides $f(U,V)$. 
		(If $(d)$ is a sharply $(q-\lambda)$-uniform direction then $S(d)\neq 0$ and this shows that $V-d$ cannot divide $f(U,V)$, so in this case the pencil with carrier $(d)$ cannot be contained in the envelope.)
		To do this, we will show 
		\[(V-d)^{\lambda-\lambda_d} \mid S(V),S_1(V), S_2(V),\ldots,S_{\lambda}(V).\]
		For $1\leq k \leq \lambda$ put $H(d)^{(k)}$ to denote the $k$-th column of $H(d)$. 
		For the integer $k$, $1 \leq k \leq \lambda-\lambda_d $\ , by \eqref{main} applied for the integers $j$ where $\lambda-k-\lambda_d\leq j \leq 2\lambda-k-\lambda_d-1$, 
		in $H(d)$ we obtain 
		that $H(d)^{(k)}$ is the linear combination of $H(d)^{(k+1)},H(d)^{(k+2)},\ldots,H(d)^{(k+\lambda_d)}$. Hence the column space of $H(d)$ is generated by the last $\lambda_d$ columns of $H(d)$, so the rank of this matrix is at most $\lambda_d$. 
		By \eqref{mainmatrix}, the column $C(d)$ is the linear combination of the columns of $H(d)$ and hence when $H(d)$ has rank at most $\lambda_d$, then so is $H_i(d)$ for each $i$.
		Then for $s>\lambda_d$ the $s\times s$ minors of $H(d)$ have determinant zero. Using Lemma 2.3 of \cite{ZsW} it follows that $d$ is a $(\lambda-\lambda_d)$-fold root of $\det H(d)=S(d)$ and of $\det H_i(d)=S_i(d)$ as we claimed.
		
		If $\lambda_d=0$ for some $(d)\in \cE_\lambda$, then $M_d$ contains the same number of copies of each element of $\F_q$ and hence $\pi_0(d)=\ldots=\pi_{2\lambda-2}(d)=0$, so $H(d)$ is the $\lambda\times \lambda$ null-matrix. From Lemma 2.3 of \cite{ZsW}, it follows that $d$ is a $\lambda$-fold root of $\det H(d)=S(d)$ and of $\det H_i(d)=S_i(d)$, thus $(V-d)^\lambda$ divides $f(U,V)$.
		
	\end{proof}
	
	The next example shows the sharpness of Theorem \ref{lambda2}.
	
	\begin{example}
		In $\AG(2,q)$, $q$ even, let $U$ denote a subset consisting of $\lambda$ points of a parabola. 
		Assume that there are at least $\lambda$  directions $(u_i)\neq (\infty)$ incident with $\lambda$ $1$-secants of $U$ (this always happens if $\lambda(\lambda-1)/2 < q - \lambda$).
		For $(d)\in \ell_{\infty}\setminus \{(\infty)\}$ denote by $N_d$ the number of $2$-secants of $U$ incident with $(d)$. 
		Then the directions of $\ell_{\infty} \setminus \{(\infty)\}$ are $(q-\lambda)$-uniform (with $0$ modulo $2$ as typical intersection number, and the renitent lines are the $1$-secants of $U$). More precisely, the direction $(d)$ is 
		sharply $(q-\lambda_d)$-uniform, where $\lambda_d=\lambda-2N_d$.
		
		Then the product of the pencils centered at $(d)$, for each $(d)$ with $N_d>0$, taking each of them with multiplicity $\lambda-\lambda_d=2N_d$, is an envelope of class 
		\[2\sum_{d} N_d=\lambda(\lambda-1).\]
		Recall that there were at least $\lambda$ sharply $(q-\lambda)$-uniform directions and so the points of $U$ are incident with at least $\lambda$ renitent lines with a sharply $(q-\lambda)$-uniform slope. Hence the renitent lines with a sharply $(q-\lambda)$-uniform direction cannot be covered by an envelope of class less than $\lambda$. But they can be covered by an envelope of class $\lambda$: the product of $\lambda$ pencils centered at the points of $U$. In total we obtain an envelope of class $\lambda(\lambda-1)+\lambda=\lambda^2$ and there is no envelope with smaller class with the given properties.
	\end{example}
	
	Consider the envelope ensured by Theorem \ref{lambda2}. The renitent lines incident with sharply $(q-\lambda)$-uniform directions of $\cE_\lambda$ are contained in a factor of degree at least $\lambda$ (since the sharply $(q-\lambda)$-uniform directions are incident with $\lambda$ of these lines, and pencils centered at such directions are not contained in the envelope) and hence from the second part of the previous theorem the following is clear. 
	
	\begin{corollary}
		\label{uccsocor}
		Suppose that the assumptions of the previous theorem holds. 
		Then there are at most $\lambda^2-\lambda$ directions incident with less than $\lambda$ renitent lines. 
		More precisely, if $\lambda_d$ denotes the number of renitent lines with slope $d$ for every $(d)\in \cE_{\lambda}$, then
		\[\sum_{(d)\in \cE_{\lambda}} (\lambda-\lambda_d)\leq \lambda^2-\lambda. \]
		It is also immediate that the number of renitent lines is at least $\lambda|\cE_\lambda|-(\lambda^2-\lambda)=\lambda(|\cE_\lambda|+1-\lambda)$.
		\qed
	\end{corollary}
	
	\section{The resultant method}
	\label{sec:4}
	
	The next result was developed in a series of papers by Sz\H{o}nyi and Weiner \cite{SzT,ZsW}, see also \cite{polybook} and the Appendix of \cite{HThesis} for the version that we cite here.
	
	\begin{result}[Sz\H onyi--Weiner Lemma]
		\label{SzWLemma}
		Let $f,g \in \F[X,Y]$ be polynomials over the arbitrary field $\F$.
		Assume that the coefficient of $X^{\deg f}$ in $f$ is not $0$ and for $y\in \F$ put $k_y = \deg \gcd(f(X,y),g(X,y))$. Then for any $y_0 \in \F$
		\begin{equation}
			\label{szwlemma}
			\sum_{y\in \F}(k_y-k_{y_0})^+ \leq (\deg f - k_{y_0})(\deg g - k_{y_0}).\ \  \footnote{Here $\alpha^+=\max\{0,\alpha\}$. Note that $g$ can be the zero polynomial as well, in that case $\deg f=k_y=k_{y_0}$ and the lemma claims the trivial $0 \leq 0$. }
		\end{equation}
	\end{result}
	
	The main ingredient of the next proofs is how we define the polynomials $f(X, Y)$ and $g(X,Y)$ in the above lemma. 
	In order to be able to detect the renitent lines at the $(q-\lambda)$-uniform directions, 
	in the definition of $g$ we introduce an auxiliary polynomial $h$ that we obtain by interpolation.
	
	This method does not reveal the underlying envelope, but it can be used to deduce quantitative information about the reintent lines.
	First we show an alternative proof of Corollary \ref{uccsocor}. Note that in this proof we do not need the $\lambda\leq (q-1)/2$ condition which was used in Theorem 		\ref{lambda2} to construct the envelope. Then we will show another quantitative result: we prove that each affine point is incident either with just a ``few" or with a ``lot" of renitent lines.

	\emph{Second proof of Corollary \ref{uccsocor}.}
	Put $|\cE_{\lambda}|=k$, $\cE_{\lambda}=\{(1:d_i:0)\}_{i=1}^k$ and $\cT=\{(a_i:b_i:1)\}^{|\cT|}_{i=1}$. 
	Denote by $m_i$ the typical intersection number corresponding to the point $(1:d_i:0)$.
	We will need the polynomial $h(Y):=\sum_{i=1}^k m_i (1 - (Y-d_i)^{q-1}) \in \F_q[Y]$. Then $h(d_i)\equiv m_i \pmod p$ for each $1\leq i \leq k$.
	%$(d)\in \cE_{\lambda}$. 
	Next define the polynomials
	\[f(X,Y):=X^q-X,\]
	\[g(X,Y):=\sum_{i}(X+a_iY-b_i)^{q-1}-|\cT|+ h(Y).\]
	Note that $g$ encodes the intersection numbers of the lines $[d:-1:x]$ with $\cT$, i.e. $g(x,d) \equiv h(d) - |[d:-1:x] \cap \cT| \pmod p$.
	For a direction $(y)\in \cE_{\lambda}$ let $\lambda_y$ denote the number of renitent lines incident with $(y)$. Then for any $(y)\in \cE_\lambda$
	it holds that 
	\[k_y:=\deg \gcd (f(X,y),g(X,y))=q-\lambda_y.\]
	Then for any fixed $(r)\in \cE_{\lambda}$, \eqref{szwlemma} gives 
	\[\sum_{(y)\in \cE_\lambda}(\lambda_r - \lambda_y) \leq \sum_{y\in \F_q}(\lambda_r - \lambda_y)^+ \leq \lambda_r (\lambda_r-1).\]
	If we choose $(r)$ so that  $\lambda_r = \lambda$, then
	\[k \lambda - \sum_{(y)\in \cE_\lambda}\lambda_y \leq \lambda (\lambda-1)\]
	and hence
	\[\sum_{(y)\in \cE_\lambda}(\lambda-\lambda_y) \leq \lambda^2-\lambda \quad \mbox{and} \quad \lambda(k+1-\lambda)\leq \sum_{(y)\in \cE_\lambda}\lambda_y.\]
	\qed	
	\begin{theorem}
		\label{resultansos}
		Take a multiset $\cT$ of $\AG(2,q)$, $q>2$, and fix an integer $\lambda>0$. 
		Let $\cF_{\lambda}$ denote the set of $(q-\lambda)$-uniform directions. 
		If $|\cF_{\lambda}|> \lambda^2+\lambda$ then for each point $R$ of the plane it holds that $R$ is incident with at most $\lambda$ or with at least  $|\cF_{\lambda}|+1-\lambda$ renitent lines. % with slope in $\cE_{\lambda}$.
	\end{theorem}
	\begin{proof}
		First we prove the following. If $\cE_{\lambda}$ is a set of at most $q$ directions which are $(q-\lambda)$-uniform and $|\cE_{\lambda}|>\lambda^2+\lambda$, then
		for each point $R$ of the plane it holds that $R$ is incident with at most $\lambda$ or with at least  $|\cE_{\lambda}|+1-\lambda$ renitent lines with slope in $\cE_{\lambda}$.

		Consider $\AG(2,q)$ embedded in $\PG(2,q)$ 	as $\{( a : b : 1) : a,b\in \GF(q) \}$ and apply a collineation of $\PG(2,q)$ which maps the line $[0 : 0 : 1]$ (the line at infinity)
		to the line  $[1 : 0 : 0]$ (the $Y$-axis) such that $(0 : 1 : 0)$ is not in the image of $\cE_{\lambda}$. 
		It may be further assumed that this collineation maps $R$ to some point $(1 : y_0 : 0)$. 
		Denote the image of $\cT$ by $\{(a_i : b_i : 1)\}_i \cup \{(1 : z_j : 0)\}_j$ and the image of $\cE_{\lambda}$ by $\{(0 : c_k : 1)\}_k$. 
		
		Let $m_k$ be the typical intersection number corresponding to the point $(0:c_k:1)$.
		We will need the polynomial $h(X)=\sum_{k=1}^{|\cE_{\lambda}|} m_k (1 - (X-c_k)^{q-1}) \in \F_q[X]$. Then $h(c_k)\equiv m_k \pmod p$ for each $k\in \{1,\ldots,|\cE_{\lambda}|\}$.
		Next define the polynomials
		\[f(X,Y):=\prod_{k=1}^{|\cE_{\lambda}|} (X-c_k),\]
		\[g(X,Y):=\sum_{i}(X+a_iY-b_i)^{q-1}+
		\sum_{j}(Y-z_j)^{q-1}-|\cT|+ h(X).\]
		As before, $g$ encodes the intersection multiplicities of the lines $[d:-1:x]$ with $\cT$, i.e. $g(x,d) \equiv h(d) - |[d:-1:x] \cap \cT| \pmod p$.
		For a point $P\in \PG(2,q)$ let $\ind P$ denote the number of renitent lines incident with $P$. Then for any $y\in \F_q$ it holds that 
		\[k_y:=\deg \gcd (f(X,y),g(X,y))=|\cE_{\lambda}|-\ind (1 : y : 0).\]
		Then \eqref{szwlemma} gives 
		\[\sum_{y\in \F_q}(\ind R - \ind (1 : y : 0))^+ \leq 
		\ind R (q-1-|\cE_{\lambda}|+\ind R).\]
		Note that 
		\[\sum_{y\in \F_q} \ind (1 : y : 0) \leq |\cE_{\lambda}| \lambda\]
		and hence
		\[q \ind R - |\cE_{\lambda}|\lambda \leq \ind R (q-1-|\cE_{\lambda}|+\ind R),\]
		\[ 0 \leq (\ind R)^2-\ind R (|\cE_{\lambda}|+1)+|\cE_{\lambda}|\lambda.\]
		For $\ind R = \lambda+1$, or $\ind R = |\cE_{\lambda}|-\lambda$ we have
		\[(\ind R)^2-\ind R (|\cE_{\lambda}|+1)+|\cE_{\lambda}|\lambda = \lambda^2+\lambda-|\cE_{\lambda}|,\]
		which is less than $0$ since $|\cE_{\lambda}|>\lambda^2+\lambda$. This proves $\ind R \leq \lambda$ or $\ind R \geq |\cE_{\lambda}|-\lambda+1$.
		
		Of course, if $|\cF_{\lambda}|\leq q$ then one can take $\cE_{\lambda}=\cF_{\lambda}$ and this proves the theorem.
		
		Now assume $q+1=|\cF_{\lambda}|>\lambda^2+\lambda$ and take an affine point $R$. We have to show that $R$ is incident with at most $\lambda$ or with at least $q+2-\lambda$ renitent lines. 
		Define $\cE_{\lambda}$ as any subset of directions of size $q$. Note that $q+1=|\cF_{\lambda}| \neq \lambda^2+\lambda+1$ because this would imply
		$q=\lambda(\lambda+1)$ and hence $\lambda=1$ and $q=2$, which we excluded. It follows that $|\cE_{\lambda}|=q > \lambda^2+\lambda$ and hence the 
		arguments above show that $R$ is incident with at most $\lambda+1$ renitent lines or with at least 
		$q+1-\lambda$ renitent lines. We have to exclude the cases when $R$ is incident with exactly $\lambda+1$ renitent lines or with exactly $q+1-\lambda$ renitent lines. 
		
		If the former case holds, then take a direction $S$ such that $RS$ is not renitent and define $\cE_{\lambda}$ as $\ell_{\infty} \setminus \{S\}$. 
		Then the renitent lines incident with $R$ have slopes in $\cE_{\lambda}$, a contradiction since there are more than $\lambda$ but less than $q+1-\lambda$ of them. 
		
		If the latter case holds, then take a direction $S$ such that $RS$ is renitent and, as before, define $\cE_{\lambda}$ as $\ell_{\infty} \setminus \{S\}$. 
		Then there are exactly $q-\lambda$ renitent lines incident with $R$ with slopes in $\cE_{\lambda}$, a contradiction since this number should be at least $q+1-\lambda$ or at most $\lambda$ (recall that $q+1 = |\cF_\lambda| > \lambda^2+\lambda$). 
	\end{proof}

	With some further efforts, the earlier methods by Sz\H onyi and Weiner \cite{SzT}, \cite{ZsW}, see also \cite{polybook}, \cite{BSzW}, also provide an alternative, but more laborious, route to reach similar results as in Section \ref{sec:3}.

	\section{\texorpdfstring{The dual statements and codewords of $\PG(2,q)$}{The dual statements and codewords of PG(2,q)}}
	\label{dualsection}
	
	In this section, we will dualize some of the earlier theorems and we will show how they apply for codewords of $\PG(2,q)$. The dual version of Theorem \ref{regular} (allowing also $|\cE|=q+1$ and not only $|\cE|\leq q$, as a consequence of Theorem \ref{gcds}, see Remark \ref{slightgen}) is the following.
	
	\begin{theorem}
		\label{regulardual}
		Fix a point $Q$ of $\PG(2,q)$. Take a multiset of lines $\cL$ in $\PG(2,q)$, none of them containing $Q$. Let $\cE_{\lambda}$ denote a non-empty set of lines through $Q$. Assume that
		\begin{enumerate}[\rm(i)]
			\item $0<\lambda \leq \min\{q-2,p-1\}$,
			\item for each line $z\in \cE_{\lambda}$, there are at most $\lambda$ points which are incident with $t_z$ modulo $p$ lines of $\cL$ (counted with multiplicities); we call these points {\em renitent},
			\item for each $z\in \cE_{\lambda}$ the remaining points of $z \setminus \{Q\}$ are incident with $m_z$ (the typical covering number of $z$) modulo $p$ lines of $\cL$, and the values $t_z-m_z$ modulo $p$ do not depend on the choice of $z$.
		\end{enumerate}
		Then the renitent points incident with lines of $\cE_{\lambda}$ are contained in an algebraic curve of degree $\lambda$.
	\end{theorem}
	
	Note that one may formulate the dual version of Theorems \ref{gcds} and \ref{lambda2} similarly.
	
	\begin{definition}
		Let $C_1(2,q)$ be the $p$-ary linear code generated by the incidence vectors of the lines of $\PG(2,q)$, $q=p^h$, $p$ prime. 
	\end{definition}
	
	The next result proves a nice geometric property on codewords of $\PG(2,q)$.
	
	\begin{result}[Blokhuis, Brouwer, Wilbrink {\cite[Proposition, p. 66]{BBW}}]
		\label{BBW}
		Let $X$ be a subset of $\PG(2, q)$ whose characteristic vector is a codeword of $C_1(2,q)$ and let $Q$ be a point not in $X$. 
		Then the points $P$ for which the line $QP$ is tangent to $X$ (i.e., $QP\cap X = \{P\}$) are all collinear.
	\end{result}
	
	The next theorem follows easily from Theorem \ref{regulardual} and it is the generalization of Result \ref{BBW}.
	
	\begin{theorem}
		\label{codewordoncurve}
		Let $X$ be a subset of $\PG(2, q)$ whose characteristic vector is a codeword $c$ of $C_1(2,q)$ and let $Q$ be a point in $\PG(2,q)$. 
		Let $\cE_{\lambda}$ denote a non-empty set of lines through $Q$, each of them incident with at most $\lambda \leq \min\{q-2,p-1\}$ points of $X\setminus\{Q\}$. 
		Then the points of $X \setminus \{Q\}$ incident with lines of $\cE_{\lambda}$ are on a curve of degree $\lambda$. 
	\end{theorem}
	
	\begin{proof}
		Put $\cE_{\lambda}=\{z_1,\ldots,z_{|\cE_{\lambda}|}\}$. 
		By definition, $c$ is the linear combination of some lines $e_j$ of $\PG(2, q)$, i.e. $c = \sum_j c_je_j$, $c_j \in \mathbb{F}_p^*$. 
		We define a multiset $\cL$ of lines of $\PG(2,q)$ in the following way: if $e_j$ is not incident with $Q$ then the weight of $e_j$ is $c_j$, all other lines have weight $0$.
		On each line $z_i$ there will be at most $\lambda$ points of $X\setminus \{Q\}$ incident with $t_{z_i}$ modulo $p$ lines of $\cL$ (counted with multiplicities) and the rest of the points will be incident with
		$m_{z_i}:=(t_{z_i}-1)$ modulo $p$ lines of $\cL$. Here $t_{z_i}=1$ and $m_{z_i}=0$ if $z_i$ is not one of the $e_j$s and $t_{z_i}=1-c_j$, $m_{z_i}=-c_j$ if $z_i=e_j$ for some $j$. Hence the result follows from Theorem \ref{regulardual}.
	\end{proof}
	
	Note that from the theorem above, it follows that the assumption $\lambda \leq p-1$ in Theorem \ref{regulardual} and so in Theorems \ref{regular} and \ref{gcds} is sharp at least in case $q=p^2$. 
	Blokhuis, Brouwer and Wilbrink (\cite{BBW}) showed that the characteristic vector of a Hermitian unital in $\PG(2, q)$ is a codeword in $C_1(2,q)$. Let $\cH$ be a Hermitian unital in $\PG(2, p^2)$ and let $Q$ be a point of $\cH$, $\cE_{\lambda}$ be the 
	set of $p^2$ lines through $Q$ intersecting the unital in $p+1$ points. If Theorem \ref{codewordoncurve} was true with $\lambda = p$, then it would yield that the points of $\cH \setminus \{ Q \}$ are on a curve $\cC$ of degree $p$. But through any point of $\cH \setminus \{ Q \}$, there pass $p^2-1$ lines intersecting $\cH \setminus \{ Q \}$ in $p+1$ points. Hence by B\'ezout's theorem these lines would be linear factors of the curve $\cC$; which is a contradiction since the degree of $\cC$ is $p < p^2-1$. 
	
	\section*{Acknowledgement}
	The authors are grateful to the anonymous referees for their valuable comments and suggestions which have certainly improved the quality of the manuscript.

	\begin{flushleft}
		Bence Csajb\'ok \\
		ELKH--ELTE Geometric and Algebraic Combinatorics Research Group\\
		ELTE E\"otv\"os Lor\'and University, Budapest, Hungary\\
		Department of Computer Science\\
		1117 Budapest, P\'azm\'any P.\ stny.\ 1/C, Hungary\\

		\medskip
		
		Current address: \\
		Dipartimento di Meccanica, Matematica e Management, \\
		Politecnico di Bari, \\
		70125 Bari, Via Orabona 4, Italy\\
		{{\tt bence.csajbok@poliba.it}}
		\medskip
		
		Peter Sziklai\\
		E\"otv\"os Lor\'and University, Budapest, Hungary\\
		Department of Computer Science\\
		1117 Budapest, P\'azm\'any P.\ stny.\ 1/C, Hungary, and\\
		R\'enyi Institute of Mathematics\\
		1053 Budapest, Re\'altanoda u. 13-15.\\
		{{\tt peter.sziklai@ttk.elte.hu}}
		
		\medskip
		
		Zsuzsa Weiner\\
		ELKH--ELTE Geometric and Algebraic Combinatorics Research Group\\
		ELTE E\"otv\"os Lor\'and University, Budapest, Hungary\\
		Department of Computer Science\\
		1117 Budapest, P\'azm\'any P.\ stny.\ 1/C, Hungary\\
		
		{{\tt zsuzsa.weiner@gmail.com}}\\
		
		and\\
		
		Prezi.com\\
		H-1065 Budapest, Nagymez\H o utca 54-56, Hungary\\
	\end{flushleft}
	

\begin{thebibliography}{pippo}
		\bibitem{survey}
		{\sc S. Ball:} Polynomials in finite geometries, in: {\it Surveys in combinatorics}, 1999, Lamb, J. D. (ed) et al, Cambridge University Press. Lond.\ Math.\ Soc.\ Lect.\ Note Ser.\ {\bf 267}, 17--35, 1999.
		
		\bibitem{arcsurvey}
		{\sc S. Ball, M. Lavrauw:} Arcs in finite projective spaces, {\it EMS Surveys in Mathematical Science} {\bf 6} (2019), 133--172.
		
		\bibitem{seminuclear}
		{\sc A. Blokhuis:} Characterization of seminuclear sets in a finite projective plane,
		{\it J.\ Geom.} {\bf 40} (1991), 15--19.
		
		\bibitem{BBW}
		{\sc A. Blokhuis, A.E. Brouwer, H. Wilbrink:} Hermitian unitals are code words, {\it Discrete Math.} {\bf 97} (1991), 63--68.
		
		\bibitem{semi}
		{\sc A. Blokhuis, T. Sz\H onyi:} Note on the structure of semiovals, {\it Discrete Math.} {\bf 106/107} (1992), 61--65.
		
		\bibitem{BSzW}
		{\sc A. Blokhuis, T. Sz\H onyi,  Zs. Weiner:} On Sets without Tangents in Galois Planes of Even Order, {\it Designs, Codes and Cryptography} {\bf 29} (2003), 91--98.
		
		\bibitem{CsW}
		{\sc B. Csajb\'ok, Zs. Weiner:} Generalizing Korchm\'aros--Mazzocca arcs, to appear in {\it Combinatorica}, {\bf 41} (2021), 601--623.
		
		\bibitem{Hilton-Milner}
		{\sc A. J. W. Hilton, E. C. Milner:} Some intersection theorems for systems of finite sets, {\it Quart. J.\ Math.\ Oxford Ser.} (2) {\bf 18} (1967), 369--384.
		
		\bibitem{HThesis}
		{\sc T. H\'eger:} Some graph theoretic aspects of finite geometries, PhD Thesis, E\"otv\"os Lor\'and University, 2013, 
		\url{http://heger.web.elte.hu//publ/HTdiss-e.pdf}
		
		\bibitem{Seidenberg} {\sc A. Seidenberg}, {\em Elements of the theory of algebraic curves}, Addison--Wesley,
		Reading, Mass., 1968.
		
		\bibitem{Simonovits}
		{\sc M. Simonovits}, Some of my Favorite Erd\H{o}s Theorems and Related Results, Theories, in: Paul Erd\H{o}s and His Mathematics, Vol. II (eds: G. Hal\'asz, L. Lov\'asz,M. Simonovits, V. T. S\'os), Springer, 2002, pp. 565--635 
		
		\bibitem{dirsziklai}{\sc P. Sziklai}, On subsets of $\mathrm{GF}(q^2)$ with $d$th power differences, {\it Discrete Math.} {\bf 208/209} (1999), 547--555.
		
		\bibitem{polybook}
		{\sc P. Sziklai:} Polynomials in Finite Geometry, Manuscript available online at \url{www.academia.edu/69422637/Polynomials_in_finite_geometry}
		
		\bibitem{SzT}
		{\sc T. Sz\H onyi:} On the Embedding of $(k,p)$-Arcs in Maximal Arcs, {\it Designs, Codes and Cryptography} {\bf 18} (1999), 235--246.
		
		\bibitem{dirszonyi}
		{\sc T. Sz\H{o}nyi:} On the number of directions determined by a set of points in an affine Galois plane,
		{\it J.\ Combin.\ Theory Ser.\ A} {\bf 74} (1996), 141--146.
		
		
		\bibitem{kmodp}
		{\sc T. Sz\H{o}nyi, Zs. Weiner:}
		Stability of $k$ mod $p$ multisets and small weight codewords of the code generated by the lines of $\PG(2,q)$,
		{\it J.\ Combin.\ Theory Ser.\ A} {\bf 157} (2018), 321--333.
		
		\bibitem{ZsW}
		{\sc Zs. Weiner:} On $(k,p^e)$-arcs in Desarguesian planes, {\it Finite Fields Appl.} {\bf 10} (2004), 390--404.
		
		
	\end{thebibliography}
\end{document}